\documentclass[12pt]{article}
\usepackage[T2A]{fontenc}
\usepackage[cp866]{inputenc}
\usepackage[english,russian]{babel}
\usepackage[tbtags]{amsmath}
\usepackage{amsfonts,amssymb,mathrsfs,amscd,amsthm}
\usepackage[lmargin=.75in,rmargin=.75in,tmargin=1.in,bmargin=1in]{geometry}
\theoremstyle{plain}

\theoremstyle{plain}

\newtheorem{theorem}{Теорема}
\newtheorem{theoremA}{Теорема~A}

\theoremstyle{definition}
\newtheorem{definition}{Определение}

\let\leq\leqslant\let\geq\geqslant

\def\supp{\operatorname{supp}}
\def\mdeg{\operatorname{deg}}

\let\myo\overline\let\myh\widehat\let\myt\widetilde
\let\leq\leqslant
\let\geq\geqslant
\def\({\left(}
\def\){\right)}

\def\FF{\mathscr F}
\def\KK{\mathscr K}
\def\mM{M}
\def\EE{E}
\def\mymu{\mu}

\def\Re{\operatorname{Re}}
\def\Im{\operatorname{Im}}
\def\mdeg{\operatorname{deg}}

\def\const{\operatorname{const}}
\def\mcap{\operatorname{cap}}
\def\RR{\mathbb R}
\def\CC{\mathbb C}

\def\RS{\mathfrak R}

\def\zz{\mathbf z}

\def\sR{\mathscr R}

\def\bad{\spaceskip=0.33emplus0.6emminus0.15em\immediate\write5{\string\bad}}
\begin{document}
\title{О сходимости аппроксимаций Чебышёва--Паде для вещественных
алгебраических функций}

\author{А.\,А.~Гончар, Е.\,А.~Рахманов, С.\,П.~Суетин}

\maketitle

\begin{abstract}
В работе анонсируются новые результаты о сходимости аппроксимаций
Чебышёва--Паде для вещественных алгебраических функций, заданных на отрезке
$[-1,1]$. Скорость сходимости аппроксимаций на отрезке и в соответствующей
``максимальной'' области мероморфности заданной функции характеризуется в
терминах некоторой теоретико-потенциальной задачи равновесия.

Библиография: 25 названий.
\end{abstract}

\markright{Аппроксимации Чебышёва--Паде}

\footnotetext[0]{Работа выполнена при поддержке Российского фонда
фундаментальных исследований (гранты \No~08-01-00317 и
\No~09-01-12160-офи-м) и Программы поддержки ведущих научных школ РФ (грант
\No~НШ-8033.2010.1).}

\subsection{}\label{s1}
Аппроксимации Паде рядов Чебышёва привлекают в последнее время особое
внимание.
Такие рациональные
аппроксимации нашли самые разнообразные применения в различных задачах
прикладной математики, теоретической физики, механики, геофизики
(см.~\cite{BaGr86}, \cite{And91}, \cite{Kniz84},~\cite{TrGu85}, \cite{Litv03},
\cite{Boy09},~\cite{Erm10}). Фактически аппроксимации Чебышёва--Паде (далее -- АПЧ) уже стали
неотъемлемой частью научных и технических расчетов, отражением чего стало
исследование оптимальных методов их вычисления (см. обзор~\cite{Litv03}) и
создание специальной программы для их нахождения в системе Maple.

В отличие от классического случая степенного ряда два стандартных способа
-- Фробениуса и Бейкера -- определения АПЧ приводят к существенно различным
результатам (см.~\cite{BaGr86}, \cite{GRS91},~\cite{GRS92}, а также п.~\ref{s2} ниже).
В~\cite{GRS91} и~\cite{GRS92} для общих ортогональных разложений были рассмотрены
оба способа построения диагональных аппроксимаций Паде и доказаны
теоремы о сходимости таких рациональных аппроксимаций для произвольной
марковской функции (см.~\eqref{3.1}). Настоящая работа является естественным
развитием
работ~\cite{GRS91} и~\cite{GRS92}: вместо марковской функции мы рассматриваем
здесь произвольную вещественную алгебраическую функцию $f$, заданную на
единичном отрезке $[-1,1]$ своим разложением в ряд Фурье--Чебышёва, и исследуем
сходимость соответствующих нелинейных и линейных аппроксимаций. Основной
результат работы -- теорема~\ref{t3} о сходимости по емкости диагональных АПЧ для
функций из указанного класса. Эта теорема является аналогом
известной теоремы Шталя~\cite{Sta86} о сходимости по емкости диагональных
аппроксимаций Паде для алгебраических функций.
Скорость сходимости АПЧ на самом отрезке $[-1,1]$ и в дополнении к
соответствующему функции $f$ ``стационарному'' компакту охарактеризована в
терминах смешанной (гриново-логарифмической) теоретико-по\-тенциальной задачи равновесия.
Отметим, что нелинейным и линейным АПЧ соответствуют различные
теоретико-по\-тенциальные задачи и, вообще говоря, различные стационарные
компакты. Эти компакты совпадают только в исключительных случаях, например, для
марковской функции $f=\myh{\sigma}$, носитель меры которой
-- отрезок вещественной прямой (см. ниже теорему~A).

\subsection{}\label{s2}
Приведем основные определения и обозначения. Пусть $\EE=[-1,1]$, $T_n(x)$ --
полиномы Чебышёва, ортонормированные на~$\EE$ с весом $(1-x^2)^{-1/2}$,
$f$ -- произвольная суммируемая на $E$ вещественная функция, заданная своим
разложением в ряд Фурье--Чебышёва:
\begin{equation}
f(x)=\sum_{k=0}^\infty c_kT_k(x),
\qquad c_k=c_k(f)=\int_{-1}^1 f(x)T_k(x)\frac{dx}{\sqrt{1-x^2}},
\quad k=0,1,\dots
\label{2.1}
\end{equation}
(другими словами заданы коэффициенты Фурье функции $f\in L_1$ по системе
$\{T_k(x)\}$). Если $f$ -- голоморфная на $\EE$ функция,
то ряд~\eqref{2.1} сходится локально равномерно внутри максимального
канонического (относительно $\EE$) эллипса голоморфности~$f$.

Фиксируем произвольную пару целых чисел $L,M\geq0$. Через $\sR(L,M)$ обозначим
класс всех рациональных функций вида $r=p/q$, где $p,q$ -- многочлены с
вещественными коэффициентами, $\mdeg{p}\leq{L}$, $\mdeg{q}\leq{M}$,
$q\not\equiv0$. Отметим, что число свободных параметров, определяющих функции
класса $\sR(L,M)$, равно $L+M+1$.

Напомним два стандартных способа определения
аппроксимаций Чебышёва--Паде функции~$f$ (или аппроксимаций Паде
ряда~\eqref{2.1}; см.~\cite[часть~2, \S~1.6]{BaGr86}). Первый способ основан на
нелинейной (относительно коэффициентов искомой рациональной функции) схеме
Бейкера, второй -- на линейной схеме Фробениуса.

Голоморфная на~$\EE$ рациональная функция $F_{L,M}$ класса $\sR(L,M)$,
разложение которой в ряд Фурье--Чебышёва имеет вид
$$
F_{L,M}(x)=c_0+c_1T_1(x)+\dots+c_{L+M}T_{L+M}+\dotsb,
$$
где $c_k=c_k(f)$, $k=0,1,\dots,L+M$, называется {\it нелинейной аппроксимацией
Паде} типа $(L,M)$ ряда~\eqref{2.1} (нелинейной АПЧ функции~$f$). Другими словами,
рациональная функция $F_{L,M}=p/q$ определяется из системы (нелинейных)
уравнений
\begin{equation}
c_k(F_{L,M})=c_k(f),\qquad k=0,1,\dots,L+M;
\label{2.3}
\end{equation}
подлежат определению из этой системы коэффициенты многочленов~$p$ и~$q$.
Система~\eqref{2.3} не всегда имеет решение и, тем самым, нелинейная
аппроксимация Чебышёва--Паде может не существовать. Так как полиномы Чебышёва
являются полиномами Фабера для отрезка $\EE$, то существование
нелинейной аппроксимации Чебышёва--Паде тесно связано с существованием
аппроксимации Паде степенного ряда
$\sum\limits_{k=0}^\infty c_kw^k$, $c_k=c_k(f)$, обладающей определенными
свойствами (см.~\cite{Sue80},~\cite{Ged81}, \cite{Sue09}).

{\it Линейной аппроксимацией Паде\/} типа $(L,M)$ ряда~\eqref{2.1}
(линейной АПЧ функции~$f$)
называется рациональная функция $\Phi_{L,M}$ класса $\sR(L,M)$, представимая
в виде $P/Q$, где $P$ и $Q$ -- произвольные многочлены ($\mdeg{P}\leq{L}$,
$\mdeg{Q}\leq{M}$, $Q\not\equiv0$), удовлетворяющие соотношениям
\begin{equation}
c_k(Qf-P)=0,\qquad k=0,1,\dots,L+M.
\label{2.2}
\end{equation}
Определяющая функцию $\Phi_{L,M}=P/Q$ система~\eqref{2.2} -- это система
линейных однородных уравнений относительно коэффициентов многочленов $P$
и~$Q$. Число уравнений системы равно $L+M+1$, число неизвестных равно $L+M+2$.
Поэтому система~\eqref{2.2} всегда имеет нетривиальное решение; легко видеть,
что этому решению соответствует $Q\not\equiv0$. Единственность такой
аппроксимации гарантировать нельзя (см., например,~\cite{Ibr02}).

Настоящая работа посвящена диагональным АПЧ; для простоты рассуждений мы
ограничимся аппроксимациями типа $(n-1,n)$.

Через $M(\EE)$ обозначим множество всех единичных (положительных
борелевских) мер, носители которых принадлежат $\EE$.
Пусть $K$ -- произвольный компакт со связным дополнением в $\myo{\CC}$ такой,
что $K\cap\EE=\varnothing$ и область
$D_K=\myo{\CC}\setminus{K}$ регулярна относительно решения задачи Дирихле,
$g_{K}(z,t)$ -- соответствующая области $D_K$ функция Грина с особенностью в
точке~$z=t\in D_K$.
Для меры $\mu\in M(\EE)$ определены логарифмический и гринов (по отношению к
компакту $K$) потенциалы:
$$
V^{\mu}(z)=\int_{-1}^1\log\frac1{|z-x|}\,d\mu(x),\qquad
G^{\mu}_{K}(z)=\int_{-1}^1g_{K}(z,x)\,d\mu(x),\qquad z\notin\EE
$$
(полагаем $g_K(z,x)\equiv0$ при $z\in K$, $x\in\EE$).
Пусть $\theta\geq0$ -- произвольное фиксированное число. Для фиксированного
компакта $K$ существует {\it единственная} мера
$\lambda(\theta)=\lambda_K(\theta)\in\mM(\EE)$, минимизирующая
функционал энергии
\begin{equation}
J_\theta(\mu;K)=\iint\(\theta\log\frac1{|x-t|}+g_{K}(x,t)\)\,d\mu(x)\,d\mu(t)
=\int\bigl(\theta V^\mu(x)+G^\mu_{K}(x)\bigr)\,d\mu(x)
\label{2.4}
\end{equation}
в классе всех мер $\mu\in\mM(\EE)$. Мера $\lambda(\theta)$ и только эта мера (в
классе $\mM(\EE)$) является {\it равновесной мерой} для смешанного
(гриново-логарифмического) потенциала
$\theta V^\mu(z)+G^\mu_K(z)$.
Другими словами, мера $\lambda(\theta)$ -- единственная мера из класса
$\mM(\EE)$, для которой имеет место соотношение равновесия
\begin{equation}
\theta V^{\lambda(\theta)}(x)+G^{\lambda(\theta)}_{K}(x)\equiv
w(\theta)=\const,\qquad x\in \EE,
\label{2.5}
\end{equation}
$w(\theta)=w_K(\theta)$ -- соответствующая {\it постоянная равновесия}; при
этом $J_\theta(\lambda;K)=w(\theta)$.

\subsection{}\label{s3}
В~\cite{GRS91} и~\cite{GRS92} была изучена сходимость нелинейных аппроксимаций
$F_n=F_{n-1,n}$
и аппроксимаций Фробениуса $\Phi_n=\Phi_{n-1,n}$ типа $(n-1,n)$ для общих
ортогональных разложений {\it марковских} функций
\begin{equation}
\myh{\sigma}(z)=\int_{F}\frac{d\sigma(x)}{z-x},\qquad z\in\myo{\CC}\setminus{F},
\label{3.1}
\end{equation}
где $F=[c,d]\subset\RR\setminus\EE$, $\sigma$ -- положительная
борелевская мера на $F$, $\sigma'=d\sigma/dx>0$ почти всюду (п.в.) на $F$.
Скорость сходимости последовательностей $F_n$ и $\Phi_n$ к функции
$f=\myh\sigma$ в области $D=\myo{\CC}\setminus{[c,d]}$ полностью
характеризуется в терминах равновесной меры
$\lambda(\theta)\in M(E)$ соответственно для $\theta=1$ и $\theta=3$ следующим
образом.

\begin{theoremA}[(см.~\cite{GRS91},~\cite{GRS92})]
Если $\sigma'>0$ п.в. на $F=[c,d]\subset\RR\setminus\EE$,
то локально равномерно в области $D\setminus\EE$
\begin{equation}
\lim_{n\to\infty}\bigl|(\myh\sigma-f_n)(z)\bigr|^{1/2n}
=\exp\bigl(-G_{F}^{\lambda(\theta)}(z)\bigr)<1,
\label{3.2}
\end{equation}
где $\theta=1$ для $f_n=F_n$ и $\theta=3$ для $f_n=\Phi_n$.
\end{theoremA}

Напомним (см.~\cite{Gon78}), что для наилучших в равномерной метрике на отрезке $\EE$
рациональных аппроксимаций $R_n=R_{n-1,n}$ функции $\myh{\sigma}$
соотношение~\eqref{3.2} справедливо с $\theta=0$.

Обозначим через $\mymu(Q)$ меру, ассоциированную с произвольным полиномом~$Q$:
$\mymu(Q)=\sum\limits_{\zeta:Q(\zeta)=0}\delta_\zeta$, где $\delta_\zeta$ --
мера Дирака с носителем в точке $\zeta$. Пусть $\myt\mu$ -- выметание меры
$\mu$ из области $\myo\CC\setminus{F}$ на~$F$. Тогда в условиях
теоремы~A для знаменателей $Q_n(z;\theta)$, $\theta=1,3,0$, соответствующих
рациональных функций $F_n,\Phi_n,R_n$ имеем:
$$
\frac1n\mymu(Q_n(\cdot;\theta))\to{\myt\lambda}(\theta),\qquad n\to\infty,
$$
где сходимость мер понимается в слабой топологии.

Отметим, что развитые в~\cite{GRS91},~\cite{GRS92} методы позволяют легко
доказать аналог теоремы~A (с заменой равномерной
сходимости~\eqref{3.2} на сходимость по емкости) и для случая, когда $F$
состоит из нескольких отрезков, а $f=\myh\sigma+r$, где $\myh\sigma$~--
марковская функция~\eqref{3.1}, $r$ -- вещественная рациональная функция,
голоморфная на~$\EE$.

\subsection{}\label{s4}
{\bad
Введем определения и обозначения, связанные с классом аналитических функций,
рассматриваемых в настоящей работе.

Компакт $K$ со связным дополнением $D_K$ в $\myo{\CC}$ будем называть {\it
допустимым} для заданной вещественной алгебраической функции~$f$, голоморфной на
отрезке~$\EE$, если
$K\cap\EE=\varnothing$ и $f$ продолжается с отрезка
$\EE$ в область $D_K$ как однозначная мероморфная функция. Множество всех
допустимых компактов для функции~$f$ обозначим через $\KK_f$.

Через $\FF(\EE)$ обозначим класс вещественных алгебраических функций~$f$,
голоморфных на $\EE$ и удовлетворяющих следующим двум условиям:

(1) существует конечное множество различных точек
$\Sigma_f=\{b_1,\dots,b_{m}\}
\subset\myo{\CC}\setminus\EE$,
$m=m(f)\geq2$, такое, что:
$\Sigma_f$ симметрично относительно вещественной прямой;
функция $f$ продолжается (с отрезка $\EE$) как многозначная аналитическая
функция в область $\myo{\CC}\setminus\Sigma_f$; каждая точка $b_j\in\Sigma_f$
является точкой ветвления функции~$f$;

(2) существует по-крайней мере один допустимый компакт $K\in\KK_f$ такой,
что $K$ симметричен относительно вещественной прямой, состоит из конечного
числа кусочно аналитических дуг и на любой открытой дуге, принадлежащей $K$,
скачок функции~$f$ отличен от тождественного нуля.

Отметим, что классу $\FF(\EE)$ принадлежат, например, следующие функции, не
являющиеся марковскими:
$$
\sqrt{(z-b)(z-\myo{b})},\quad
\root3\of{(z-b)(z-\myo{b})(z-a)},
$$
где $\Im{b}>0$, $a\in\RR\setminus\EE$ и выбрана надлежащая ветвь алгебраической
функции.}

В дальнейшем функция $f\in\FF(E)$ предполагается фиксированной.

\subsection{}\label{s5}
Хорошо известно (см.~\cite{Sta86},~\cite{GoRa87}), что при
доказательстве сходимости диагональных аппроксимаций Паде для алгебраических
функций ключевую роль играет существование допустимого (для заданной функции)
компакта, обладающего так называемым свойством {\it стационарности} (или
$S$-{\it свойством}). Это понятие оказывается тесно связанным с соответствующей
теоретико-по\-тенциальной задачей равновесия.
Приведем определение $S$-свойства допустимого компакта, соответствующего
рассматриваемой задаче равновесия~\eqref{2.5}.

\begin{definition}
Пусть параметр $\theta\geq0$. Будем говорить, что
(не разбивающий плоскость и состоящий из конечного числа кусочно-аналитических
дуг) допустимый компакт $F=F(\theta)\in\KK_f$ обладает свойством {\it симметрии}
(или $S$-{\it свойством}), если:
\begin{equation}
\frac{\partial G^{\lambda}_{F}}{\partial n_{+}}(\zeta)
=\frac{\partial G^{\lambda}_{F}}{\partial n_{-}}(\zeta),
\qquad\zeta\in F^0,
\label{5.0}
\end{equation}
где $\lambda=\lambda_F(\theta)$ -- соответствующая равновесная мера,
$F^0$ -- объединение всех открытых дуг, принадлежащих компакту $F$,
$\partial/\partial n_{\pm}$ -- нормальные производные, взятые с противоположных
сторон $F^0$.
\end{definition}

Пусть параметр $\theta=1$, $\lambda=\lambda_K(1)\in M(E)$ --
равновесная мера, соответствующая произвольному компакту $K\in\KK_f$,
$w=w_K(1)$ -- соответствующая постоянная равновесия (см.~\eqref{2.5}):
$V^{\lambda}(x)+G^{\lambda}_{K}(x)\equiv w$, $x\in \EE$;
при этом $J(\lambda;1)=\min\limits_{\mu\in\mM(\EE)}J(\mu;1)=w$.
Справедлива следующая
\begin{theorem}\label{t1}
Если функция $f\in\FF(\EE)$, то существует единственный компакт
$F=F(1)\in\KK_f$ такой,
что
\begin{equation}
J(\lambda_{F};1)=\max_{K\in\KK_f}J(\lambda_K;1).
\label{5.1}
\end{equation}
Стационарный компакт $F$ состоит из конечного числа кусочно-аналитических дуг, не
разбивает плоскость и обладает $S$-свойством~\eqref{5.0},
где $\lambda=\lambda(1)$ -- соответствующая равновесная мера.
\end{theorem}

Теорема~\ref{t1} доказывается в два этапа в соответствии со следующей схемой.
Сначала с помощью геометрических соображений, основанных
на замене при $\theta=1$ функционала энергии~\eqref{2.4} меры $\lambda\in\mM(\EE)$ на
обобщенный (по отношению к допустимому компакту $K\in\KK_f$)
трансфинитный диаметр~$\EE$ доказывается, что максимум в правой
части~\eqref{5.1} достаточно искать среди тех допустимых компактов, которые
лежат вне максимального канонического эллипса голоморфности~$f$. Такое
семейство компактно в хаусдорфовой метрике, поэтому существует допустимый
компакт~$F$, удовлетворяющий соотношению~\eqref{5.1}. Затем с помощью
вариационного метода аналогично~\cite{PeRa94} устанавливается, что этот экстремальный
компакт~$F$ является замыканием критических траекторий некоторого
квадратичного дифференциала. Отсюда уже вытекает $S$-свойство~\eqref{5.0}.

Пусть теперь параметр $\theta=3$, $\lambda=\lambda_K(3)\in M(E)$ --
равновесная мера, соответствующая произвольному компакту $K\in\KK_f$,
$w=w_K(3)$ -- соответствующая постоянная равновесия:
$3V^{\lambda}(x)+G^{\lambda}_{K}(x)\equiv w$, $x\in \EE$;
при этом $J(\lambda;3)=\min\limits_{\mu\in\mM(\EE)}J(\mu;3)=w$.
Пусть $\rho(f)$ -- индекс максимального канонического эллипса, в который функция
$f$ продолжается с отрезка $E$ как голоморфная функция.
Справедлива следующая

\begin{theorem}\label{t2}
Если функция $f\in\FF(\EE)$ и $\rho(f)>\sqrt{2}$, то существует единственный
компакт $F=F(3)\in\KK_f$ такой, что
\begin{equation}
J(\lambda_{F};3)=\max_{K\in\KK_f}J(\lambda_K;3).
\label{6.1}
\end{equation}
Стационарный
компакт $F$ состоит из конечного числа кусочно-аналитических дуг, не
разбивает плоскость и обладает $S$-свойством~\eqref{5.0},
где $\lambda=\lambda(3)$ -- соответствующая равновесная мера.
\end{theorem}

Отметим, что параметрам $\theta=1$ и $\theta=3$ соответствуют существенно
разные векторные теоретико-потенциальные задачи равновесия (см. ниже
п.~\ref{s8}); ограничение $\rho(f)>\sqrt2$ связано именно с этим различием.

Справедлива следующая
\begin{theorem}\label{t3}
Пусть $f\in\FF(\EE)$.
Тогда для любого компакта $K\subset\myo{\CC}\setminus(\EE\cup{F})$
\begin{equation}
\bigl|(f-f_n)(z)\bigr|^{1/2n}\overset{\mcap}{\longrightarrow}
\exp\bigl(-G_{F}^{\lambda(\theta)}(z)\bigr)<1,
\qquad z\in K,
\label{5.3}
\end{equation}
где $\theta=1$, $F=F(1)$ для $f_n=F_n$ и $\theta=3$, $F=F(3)$ для $f_n=\Phi_n$.
\end{theorem}

Отметим, что для $\theta=1$ утверждение теоремы~\ref{t3} имеет место, вообще говоря,
по некоторой подпоследовательности, для которой
существуют нелинейные диагональные АПЧ
(см.~\cite{Sue80}, \cite{Sue09},~\cite{Sta96},~\cite{Sue00}).

Таким образом в условиях теоремы~\ref{t3}
$f_n\overset{\mcap}\longrightarrow f$ на компактных подмножествах
области $D=\myo{\CC}\setminus{F}$ и
$$
\frac1n\mymu(Q_n(\cdot,\theta)\to\myt{\lambda}(\theta),\qquad n\to\infty,
$$
где $\myt{\lambda}(\theta)$ -- выметание меры
$\lambda=\lambda(\theta)\in\mM(\EE)$ из области $D$ на $\partial D=F$.

Теорема~\ref{t3} доказывается по общей схеме, предложенной в~\cite{Sta86},
\cite{GoRa87} и основанной на $S$-свойстве~\eqref{5.0} соответствующего
стационарного компакта.
Отметим, что из~\eqref{5.3} вытекает, что каждый полюс~$f$ в~$D$ притягивает
по-крайней мере столько полюсов $f_n$, какова его кратность.

Таким образом, для функции~$f$ класса $\FF(\EE)$ компакты $F(1)$ и $F(3)$ в
случае соответственно нелинейных и линейных аппроксимаций Чебышёва--Паде играют
роль отрезка $[c,d]$, соответствующего марковской функции $\myh\sigma$,
$\supp{\sigma}=[c,d]$ (см.~\eqref{3.1}).

Отметим, что при $\theta=1$ стационарный компакт $F(1)$ является
образом компакта Шталя для степенного ряда $\sum\limits_{k=0}^\infty c_kw^k$,
$c_k=c_k(f)$, при отображении, задаваемом функцией Жуковского $z=(w+w^{-1})/2$.

\subsection{}\label{s6}
Пусть параметр $\theta=1$. Если соответствующий стационарный компакт $F=F(1)$
состоит из $s=s(1)$ непересекающихся аналитических
дуг, попарно соединяющих точки ветвления $b'_1,\dots,b'_{2s}\in\Sigma_f$
функции~$f$, то он допускает наглядное описание в терминах, связанных с
{\it четырехлистной} римановой поверхностью рода $g=s-1$.

Построим сначала двулистную риманову поверхность
$\RS=\RS^{(1)}\cup\RS^{(2)}$ следующим образом. Возьмем два экземпляра
римановой сферы $\myo{\CC}$, разрезанных по отрезку $\EE=[-1,1]$, и переклеим
по разрезу. Полученная двулистная риманова поверхность
$\RS=\RS^{(1)}\cup\RS^{(2)}$ эквивалентна римановой сфере. Определим на $\RS$
функцию $u(\zz)$, $\zz\in\RS$, следующим образом:
$u(z^{(1)})=G^{\lambda_{F}}_{F}(z)$, $u(z^{(2)})=w_F-V^{\lambda_{F}}(z)$.
Непосредственно из условия равновесия~\eqref{2.5} вытекает, что $u$ --
гармоническая функция на $\RS\setminus(F^{(1)}\setminus\{\infty^{(2)}\})$.
Кроме того, $u\equiv0$ на компакте $F^{(1)}$ и
$u(z^{(2)})=\log|z|+w_F+o(1)$ при $z^{(2)}\to\infty^{(2)}$.
Возьмем теперь два экземпляра $\RS$, разрезанных по компакту $F^{(1)}$, и
переклеим их между собой по
соответствующим разрезам. Получим четырехлистную риманову поверхность $\RS_1$
рода $g=s-1$. Так как $u\equiv0$ на $F^{(1)}$, то $u$
гармонически продолжается с одного экземпляра $\RS$ на другой с переменой
знака. Продолженная функция гармонична на $\RS_3$ всюду кроме точек
$\zz=\infty^{(2)},\infty^{(3)}$, где она имеет
логарифмические особенности:
$\log{|z|}$ при $\zz\to\infty^{(2)}$ и
$-\log{|z|}$ при $\zz\to\infty^{(3)}$.
Следовательно, $u(\zz)=\Re\Omega(\zz)$,
$\Omega(\zz)=\displaystyle\int_{b'_1}^{\zz}d\Omega$, где $d\Omega(\zz)$ --
(единственный) абелев
дифференциал  на $\RS_1$ с чисто мнимыми периодами и особенностью вида
$1/z$ в точке $\zz=\infty^{(2)}$ и вида
$-1/z$ в точке $\zz=\infty^{(3)}$. Компакт $F$
соответствует нулевой линии уровня
функции $\Re\Omega(\zz)$: $F=\{z\in\myo{\CC}:\Re\Omega(\zz)=0\}\setminus\EE$.

Отметим, что $u(\zz)=g_{F^{(1)}}(\zz,\infty^{(2)})$ --
функция Грина для области $\RS\setminus{F^{(1)}}$ с особенностью в точке
$\zz=\infty^{(2)}$, а $w_F=\gamma^{(2)}$ -- постоянная Робена для этой
функции Грина. Тем самым задача о максимуме
постоянной $w_F$ соответствует задаче о минимуме $e^{-\gamma}$, т.е. минимуме
логарифмической емкости.

\subsection{}\label{s7}
Если для параметра $\theta=3$ стационарный компакт $F=F(3)$ состоит из
$s=s(3)$ непересекающихся аналитических дуг, попарно соединяющих точки
ветвления $b''_1,\dots,b''_{2s}\in\Sigma_f$ функции~$f$, то он допускает
наглядное описание в терминах, связанных с {\it шестилистной} римановой
поверхностью рода $g=s-1$.

Построим сначала трехлистную риманову поверхность
$\RS=\RS^{(1)}\cup\RS^{(2)}\cup\RS^{(3)}$ следующим образом. Возьмем три экземпляра римановой
сферы $\myo{\CC}$. На первом проведем разрез по отрезку $\EE=[-1,1]$, на втором
-- по $\EE$ и дугам $\ell_j$, $j=1,\dots,s$, составляющим $F$, на третьем --
только по дугам $\ell_j$.
Три полученных экземпляра римановой сферы переклеиваются друг с
другом следующим образом. Второй подклеивается к первому по разрезу,
соответствующему отрезку~$\EE$, третий -- ко второму по разрезам, соответствующим
дугам $\ell_j$. Полученная трехлистная риманова поверхность
$\RS=\RS^{(1)}\cup\RS^{(2)}\cup\RS^{(3)}$ эквивалентна римановой сфере
с $g=s-1$ ручками.
Определим на $\RS$ функцию $u(\zz)$, $\zz\in\RS$, следующим образом:
$u(z^{(1)})=2G^{\lambda_{F}}_{F}(z)$,
$u(z^{(2)})=G^{\lambda_{F}}_{F}(z)+w_{F}-3V^{\lambda_{F}}(z)$,
$u(z^{(3)})=-G^{\lambda_{F}}_{F}(z)+w_{F}-3V^{\lambda_{F}}(z)$.
Непосредственно из условия равновесия~\eqref{2.5} и $S$-свойства~\eqref{5.0}
компакта
$F$ вытекает, что $u$ -- гармоническая функция на
$\RS\setminus(F^{(1)}\cup\EE^{(3)})\setminus\{\infty^{(2)},\infty^{(3)}\}$.
Кроме того, $u\equiv0$ на компакте $K=F^{(1)}\cup\EE^{(3)}$ и
$u(z^{(2)})=3\log|z|+c_2+o(1)$ при $z^{(2)}\to\infty^{(2)}$,
$u(z^{(3)})=3\log|z|+c_3+o(1)$ при $z^{(3)}\to\infty^{(3)}$.
Отсюда вытекает, что $u(\zz)\equiv3W(\zz)$, где
$W(\zz)=g_K(\zz,\infty^{(2)})+g_K(\zz,\infty^{(3)})$,
$g_K(\zz,\,\cdot\,)$ -- функция Грина для области $\RS\setminus{K}$ с
особенностью в соответствующей точке. Возьмем теперь два экземпляра $\RS$,
разрезанных
по отрезку $\EE^{(3)}$ и компакту $F^{(1)}$, и переклеим их между собой по
соответствующим разрезам. Получим шестилистную риманову поверхность $\RS_3$
рода $g=s-1$. Так как $u\equiv0$ на $F^{(1)}\cup\EE^{(3)}$, то $u$
гармонически продолжается с одного экземпляра $\RS$ на другой с заменой знака.
Продолженная функция гармонична на $\RS_1$ всюду кроме точек
$\zz=\infty^{(2)},\infty^{(3)},\infty^{(4)},\infty^{(5)}$, где она имеет
логарифмические особенности:
$3\log{|z|}$ при $\zz\to\infty^{(2)},\infty^{(3)}$ и
$-3\log{|z|}$ при $\zz\to\infty^{(4)},\infty^{(5)}$.
Следовательно, $u(\zz)=\Re\Omega(\zz)$,
$\Omega(\zz)=\displaystyle\int_{b''_1}^{\zz}d\Omega$, где $d\Omega(\zz)$ --
(единственный) абелев
дифференциал  на $\RS_1$ с чисто мнимыми периодами и особенностями вида
$1/z$ в точках $\zz=\infty^{(2)},\infty^{(3)}$ и вида
$-1/z$ в точках $\zz=\infty^{(4)},\infty^{(5)}$. Компакт $F$
соответствует нулевой линии уровня
функции $\Re\Omega(\zz)$: $F=\{z\in\myo{\CC}:\Re\Omega(\zz)=0\}\setminus\EE$.

Отметим, что для параметра $\theta=3$ постоянная
$w_F=\gamma(\infty^{(2)})+\gamma(\infty^{(3)})+2g_D(\infty^{(2)},\infty^{(3)})$,
где $\gamma(\infty^{(2)})$, $\gamma(\infty^{(3)})$ -- постоянные Робена в
точках $\infty^{(2)}$ и $\infty^{(3)}$ для функции Грина области
$D=\RS\setminus({F^{(1)}}\cup{E^{(3)}})$.

Для произвольной области $D\subset\CC$ величина
$\gamma(a_1)+\gamma(a_2)+2g_D(a_1,a_2)$ равна {\it приведенному модулю} области $D$
относительно двух точек $a_1,a_2$. Это понятие (для произвольного набора
точек) было введено В.\,Н.~Дубининым в~\cite{Dub94a},~\cite{Dub94b}
геометрическим образом
через модуль плоского конденсатора. Тем самым экстремальная константа $w(F)$
(т.е. максимум минимума энергии) соответствует максимальному {\it приведенному
модулю} (относительно двух точек).

\subsection{}\label{s8}
Непосредственно из результатов
работ~\cite{Sta86} и~\cite{GoRa87} вытекает, что для {\it наилучших} в равномерной
метрике на отрезке $[-1,1]$ рациональных аппроксимаций $f_n=R_n$ функции
$f\in\FF(\EE)$ также справедливо
соотношение~\eqref{5.3}, где $\theta=0$, $F=F(0)$ -- стационарный
компакт, соответствующий задаче равновесия~\eqref{2.5} с
$\theta=0$ и обладающий $S$-свойством~\eqref{5.0}, $\lambda=\lambda_F(0)$ --
соответствующая равновесная мера.
В этом случае функция $u(z^{(1)})=w_F-G^\lambda(z)$ продолжается через разрез,
проведенный по отрезку $E$, на второй лист римановой поверхности
$\RS=\RS^{(1)}\cup\RS^{(2)}$ с переменой знака. Отсюда уже легко вытекает,
что задача о максимуме постоянной $w_F$ эквивалентна задаче о минимуме
емкости конденсатора $(F^{(1)},F^{(2)})$.

Таким образом, все три стационарных (в заданном классе $\KK_f$)
компакта $F(1),F(3),F(0)$ обладают $S$-свойством~\eqref{5.0}. Это свойство, тем
самым, носит вполне универсальный характер.

Отметим в заключение, что параметрам $\theta=0$, $\theta=1$ и $\theta=3$
соответствуют существенно разные векторные (размера $2\times2$)
теоретико-потенциальные задачи равновесия. При $\theta=0$ матрица
взаимодействия имеет вид
$A_0=\begin{pmatrix} 1&-1\\-1&1\end{pmatrix}$, при $\theta=1$ матрица
$A_1=\begin{pmatrix} 2&-1\\-1&1\end{pmatrix}$, при $\theta=3$ матрица
$A_3=\begin{pmatrix} 4&-1\\-1&1\end{pmatrix}$.

Доказательствам теорем~\ref{t1},~\ref{t2},~\ref{t3} предполагается посвятить
работы~\cite{GRS11a} и~\cite{GRS11b}.

\end{document}